\documentclass[12pt,leqno,twoside]{article}

\usepackage{fancyhdr}
\usepackage{titlesec}
\usepackage{microtype}
\usepackage{amsmath}
\usepackage{amsthm}
\usepackage{amsfonts}
\usepackage{amssymb}
\usepackage{graphicx}
\usepackage{color}
\usepackage{enumerate}
\usepackage[utf8]{inputenc}
\usepackage{hyperref}

\DeclareMathOperator{\tb}{tb}
\DeclareMathOperator{\rot}{rot}
\DeclareMathOperator{\selfl}{sl}
\DeclareMathOperator{\PD}{PD}
\DeclareMathOperator{\lk}{lk}

\newcommand{\Z}{\mathbb{Z}}
\newcommand{\Q}{\mathbb{Q}}
\newcommand{\N}{\mathbb{N}}

\newtheoremstyle{thm}{.5\baselineskip}{.5\baselineskip}{\itshape}{5mm}
{\scshape}{.}{2mm}{}

\theoremstyle{thm}
\newtheorem{theorem}{Theorem}[section]
\newtheorem{corollary}[theorem]{Corollary}
\newtheorem{lemma}[theorem]{Lemma}

\newtheorem{definition}[theorem]{Definition}

\newtheoremstyle{rem}{.5\baselineskip}{.5\baselineskip}{}{5mm}
{\scshape}{.}{2mm}{}

\theoremstyle{rem}
\newtheorem{remark}[theorem]{Remark}
\newtheorem{example}[theorem]{Example}

\renewcommand{\abstractname}{{\bf Abstract.\ }}
\renewenvironment{abstract}{\footnotesize\abstractname}

\makeatletter

\renewcommand{\section}{\@startsection{section}{1}{0pt}{\baselineskip} 
{.5\baselineskip}{\centering\normalfont\bfseries}}

\renewcommand{\subsection}{\@startsection{subsection}{1}{7mm}{\baselineskip}
{-3mm}{\bfseries}}

\renewcommand*{\@seccntformat}[1]{%
  \csname the#1\endcsname.\ }

\makeatother

\begin{document}

\pagestyle{fancy}

\thispagestyle{empty}

\begin{center}
{\bf\large COMPUTING ROTATION AND SELF-LINKING NUMBERS IN CONTACT SURGERY DIAGRAMS
}

\bigskip

{\begin{NoHyper}\small SEBASTIAN DURST$^1$ and MARC KEGEL$^1$
\footnote{\parindent=6mm{\it Key words and phrases:} rotation number, Legendrian knots,
transverse knots, contact surgery diagrams. 

}\end{NoHyper}}

\medskip

{\scriptsize $^1$Mathematisches Institut, Universit\"at zu K\"oln, Weyertal 86--90, 50931 K\"oln, Germany \\[-2mm]
e-mail: sdurst@math.uni-koeln.de, mkegel@math.uni-koeln.de}

\end{center}

\bigskip

\begin{quotation}
\begin{abstract}
We give an explicit formula to compute the rotation number of a nullhomologous Legendrian knot in contact ($1/n$)-surgery diagrams along Legendrian links and obtain a corresponding result for the self-linking number of transverse knots. Moreover, we extend the formula by Ding--Geiges--Stipsicz for computing the $d_3$-invariant to ($1/n$)-surgeries.
\end{abstract}
\end{quotation}

\section{Introduction}

A lot of the geometry of a $3$-dimensional contact manifold is encoded in its \emph{Legendrian knots}, i.e.\ smooth knots tangent to the contact structure, and in its \emph{transverse knots}, i.e.\ smooth knots transverse to the contact structure. Therefore a main topic in $3$-dimensional contact geometry is the study of these knots. In particular, it is a challenge to distinguish knots within these classes. For nullhomologous knots this is mostly done by the so-called \emph{classical invariants}, the \emph{Thurston--Bennequin invariant} and the \emph{rotation number} for Legendrian knots and the \emph{self-linking number} for transverse knots. In the unique tight contact structure of the $3$-sphere there are easy formulas to compute the classical invariants from a \emph{front projection} of the knot. For this and other basic notions in contact geometry we refer the reader to \cite{Geiges2008}.

A natural extension is to consider Legendrian or transverse knots in \emph{contact surgery diagrams} along Legendrian links and to compute their classical invariants in the surgered manifold. Starting with the work of Lisca, Ozsv{\'a}th, Stipsicz and Szab{\'o} \cite[Lemma 6.6]{MR2557137} several results were obtained in that setting by Geiges and Onaran \cite[Lemma 2]{MR3338830}, Conway \cite[Lemma 6.4]{Conway2014} and Kegel \cite[Section 8]{Kegel2016}. 

In Theorem~\ref{thm:rot_surgery} we combine the results mentioned above to obtain a condition when a Legendrian knot is nullhomologous in the surgered manifold and give a formula computing its rotation number. In Section~\ref{section:selfl} we explain how to represent a transverse knot in a contact surgery diagram along Legendrian knots and then compute its self-linking number in the surgered contact manifold. Finally, in Section~\ref{section:rational} we extend these results to rationally nullhomologous knots. On the way we present plenty of examples on how to use these formulas.

A further closely related topic is the computation of the $d_3$-invariant of the resulting
contact structure in a surgery diagram. By translating contact $(1/n)$-surgeries into
($\pm 1$)-surgeries, we generalise the formula by Ding--Geiges--Stipsicz \cite{MR2056760} in Section~\ref{section:d3}.
\section{The rotation number in surgery diagrams}
\label{section:rotation}

Let $L = L_1 \sqcup \ldots \sqcup L_k \subset S^3$ be an oriented link in $S^3$ and let $M = S^3_L(r)$ be the manifold obtained by surgery along $L$ with coefficients ${p_i}/{q_i}$ (for basic notions of Dehn surgery see \cite{Rolfsen2003}). We denote the corresponding surgery slopes $r_i = p_i \mu_i + q_i \lambda_i \in H_1 (\partial \nu L_i)$, where $\mu_i$ is represented by a positive meridian of $L_i$ and $\lambda_i$ is the Seifert longitude of $L_i$.
If no coefficient group is specified, homology groups are understood to be over the integers.
Let $L_0 \subset S^3\setminus L$ be an oriented knot in the complement of $L$.

Define $l_{ij} := \lk (L_i, L_j)$ and let $\mathbf{l}$ be the vector with components $l_i = l_{0i}$ and $Q$ the generalised linking matrix:
$$
Q = \begin{pmatrix} p_1 & q_2 l_{12} & \cdots & q_n l_{1k} \\
q_1 l_{21} & p_2 & & \\
\vdots & & \ddots & \\
q_1 l_{k1} & & & p_k \end{pmatrix}.
$$
The knot $L_0$ is nullhomologous in $M$ if and only if there is an integral solution $\mathbf{a}$ of the equation $\mathbf{l} = Q\mathbf{a}$ (see \cite{Kegel2016}).

\begin{definition}
Let $K \subset (M, \xi)$ be a nullhomologous oriented Legendrian knot and $\Sigma$ a Seifert surface for $K$.
The \emph{rotation number} of $K$ with respect to the Seifert surface $\Sigma$ is equal to
$$
\rot(K,\Sigma) = \langle \textrm{e}(\xi, K), [\Sigma] \rangle = \PD(\textrm{e}(\xi, K)) \bullet [\Sigma],
$$
where $\textrm{e}(\xi, K)$ is the relative Euler class of the contact structure
$\xi$ relative to the trivialisation given by a positive tangent vector field
along the knot $K$, and $[\Sigma]$ the relative homology class represented by
the surface $\Sigma$.
\end{definition}

This definition of the rotation number is useful for calculations (see also \cite{Ozsvath2015}).
For an alternative equivalent definition see \cite{Geiges2008}. Clearly, the rotation number does only depend on the class of the chosen Seifert surface, not on the particular choice of surface itself. Note also that the rotation number is independent of the class of the Seifert surface if the Euler class $\textrm{e}(\xi)$ of the contact structure vanishes (see Proposition 3.5.15 in \cite{Geiges2008}).

\begin{theorem}
\label{thm:rot_surgery}
Let $L = L_1 \sqcup \ldots \sqcup L_k$ be an oriented Legendrian link in $(S^3, \xi_\textrm{st})$ and $L_0$ an oriented Legendrian knot in its complement. Let $(M, \xi)$ be the contact manifold obtained from $S^3$ by contact $({1}/{n_i})$-surgeries ($n_i \in \Z$) along $L$.
Then $L_0$ is nullhomologous in $M$ if and only if there is an integral vector $\mathbf{a}$ solving $\mathbf{l} = Q\mathbf{a}$ as above, in which case its rotation number in $(M, \xi)$ with respect to a special Seifert class $\widehat{\Sigma}$ constructed in the proof is equal to
$$
\rot_M (L_0, \widehat{\Sigma}) = \rot_{S^3} (L_0) - \sum_{i=1}^k{a_i n_i \rot_{S^3} (L_i)}.
$$
\end{theorem}

The proof proceeds in two steps. First, following \cite{Conway2014}, we construct the class of a Seifert surface for $L_0$ in $M$. We then use the description of the rotation number in terms of relative Euler classes to compute $\rot$.

\begin{remark}
\label{coefficients}\hfill
\begin{enumerate}
\item Notice that the matrix $Q$ is formed using the topological surgery coefficients $p_i / q_i$, not the contact surgery coefficients.
The topological surgery coefficient equals the sum of the contact surgery coefficient and the Thurston--Bennequin invariant of the surgery knot. Therefore, we always have $q_i = n_i$.
	\item
	Observe that for any contact surgery coefficient $r \neq 0$ there exists a tight contact structure on the glued in solid torus compatible with the surgery. This tight contact structure on the solid torus is unique if and only if the surgery coefficient is of the form ${1}/{n}$ for $n \in \Z$.
	Therefore, contact $({1}/{n})$-surgery is well-defined (see \cite[Proposition~7]{MR1823497}).

	For a general contact $r$-surgery, there is an algorithm transforming the surgery
	into contact $({1}/{n})$-surgeries. The procedure is not unique, however, the algorithm provides all possible choices of contact structures that are tight on the surgery torus (cf.\ \cite{MR2055048, MR2056760}).
	In contrast to the Thurston--Bennequin invariant, the rotation number in the surgered manifold does indeed depend on the choice of contact structure on the surgery tori, cf.\ Example~\ref{ex:non-unique}.
	\item In \cite{MR2055048} it is shown that one can get any contact $3$-manifold by a sequence of contact $(1/n)$-surgeries starting from the standard tight $3$-sphere. Moreover, it is easy to show that any Legendrian knot in the resulting contact manifold can be represented by a Legendrian knot in the complement of the surgery link. 
\end{enumerate}
\end{remark}

\begin{proof}
Assume that $L_0$ is nullhomologous in $M$ and fix Seifert surfaces $\Sigma_0, \ldots, \Sigma_k$ for $L_0, \ldots, L_k$ in $S^3$, such that intersections of surfaces and link components are transverse.
Our aim is to use these surfaces to construct the class of a Seifert surface for $L_0$ in the surgered manifold $M$.
By abuse of notation, we will identify $\Sigma_i$ with its class in $H_2 (S^3 \setminus \nu L_i, \partial \nu L_i)$ and will denote the class in $H_2 (S^3 \setminus (L_0 \sqcup \nu L), \partial L_0 \sqcup \partial\nu L)$ induced by restriction again by $\Sigma_i$.

The idea is to construct a class of the form
$$
\Sigma = \Sigma_0 + \sum_{i=1}^{k}{k_i \Sigma_i}
$$
such that its image under the boundary homomorphism $\partial$ in the long exact sequence of the pair $(S^3 \setminus (L_0 \sqcup \nu L), \partial L_0 \sqcup \partial\nu L)$ is a linear combination of the surgery slopes $r_i$ and a longitude of $L_0$, i.e.\ we want
$$
\partial \Sigma = t\mu_0 + \lambda_0 + \sum_{i=1}^{k}{m_i r_i} = t\mu_0 + \lambda_0 + \sum_{i=1}^{k}{m_i ( p_i\mu_i + q_i\lambda_i )}.
$$
So our aim is to solve this equation for $k$ and $m_i$. We will first describe the boundary homomorphism $\partial$ in more detail and then compare coefficients.
The surgery slopes $r_i$ bound discs in the surgered manifold $M$, so $\Sigma$ can be extended to give rise to a class in $H_2 (M\setminus \nu L_0, \partial \nu L_0)$, which we denote by $\widehat{\Sigma}$.
Geometrically, the boundary homomorphism sends $\Sigma$ to its intersection with the boundary of the link complement. So we have:
$$
\partial \colon\thinspace \Sigma_j \longmapsto \lambda_j - \sum_{i\neq j}{l_{ij} \mu_i},
$$
and thus
\begin{align*}
\partial \colon\thinspace \Sigma  \longmapsto & \sum_{i=0}^{k}{k_i \lambda_i} - \sum_{j=0}^{k}\sum_{i\neq j}{k_j l_{ij} \mu_i} \\
 = & -\sum_{j=1}^{k}{k_j l_{0j}} \mu_0 + \lambda_0 + \sum_{i=1}^{k}{k_i \lambda_i} - \sum_{i=1}^k{l_{0i} \mu_i} - \sum_{i=1}^{k}\sum_{j\neq i}{k_j l_{ij} \mu_i},
\end{align*}
where we set $k_0 = 1$. Note that the minus sign stems from the induced boundary orientation of $\Sigma$ (see Figure \ref{fig:orientations}).
Setting $k_i = -a_i q_i$ and using that $L$ is nullhomologous, we obtain
\begin{align*}
\partial \colon\thinspace \Sigma_j \longmapsto & \sum_{j=1}^{k}{a_j q_j l_{0j}} \mu_0 + \lambda_0 - \sum_{i=1}^{k}{a_i q_i \lambda_i} - \sum_{i=1}^k{l_{0i} \mu_i} + \sum_{i=1}^{k}\sum_{j\neq i}{a_j q_j l_{ij} \mu_i} \\
= & \sum_{j=1}^{k}{a_j q_j l_{j}} \mu_0 + \lambda_0 - \sum_{i=1}^{k}{a_i q_i \lambda_i} - \sum_{i=1}^k{l_{i} \mu_i} + \sum_{i=1}^{k}\sum_{j\neq i}{a_j Q_{ij} \mu_i} \\
= & \sum_{j=1}^{k}{a_j q_j l_{j}} \mu_0 + \lambda_0 - \sum_{i=1}^{k}{a_i q_i \lambda_i} - \sum_{i=1}^k{l_{i} \mu_i} + \sum_{i=1}^{k}{(l_i - a_i  p_i )\mu_i} \\
= & \sum_{j=1}^{k}{a_j q_j l_{j}} \mu_0 + \lambda_0 - \sum_{i=1}^{k}{a_i q_i \lambda_i} - \sum_{i=1}^{k}{(a_i  p_i )\mu_i},
\end{align*}
which is of the desired form with\ $m_i = -a_i$ and $t = \sum_{j=1}^{k}{a_j q_j l_{j}}$.

\begin{figure}[htbp] 
\centering
\def\svgwidth{0,88\columnwidth}
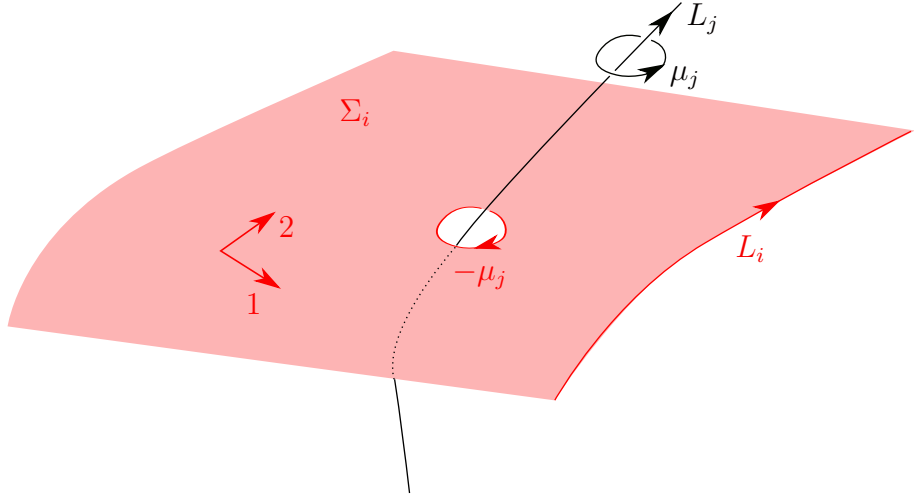
\caption{Orientation of the meridian induced by the intersection}
\label{fig:orientations}
\end{figure}	

\begin{remark}
Observe that we can also directly obtain an embedded surface representing the capped-off class $\widehat{\Sigma}$ by resolving self-intersections in $\Sigma$.
In particular, $t$ is the negative change of the Thurston--Bennequin number of $L_0$ in the surgery (cf.\ \cite{Conway2014}, \cite{Kegel2016}), i.e. we get
$$
\tb_{M} (L_0) = \tb_{S^3} (L_0) - \sum_{j=1}^{k}{a_j n_j l_{j}}.
$$
\end{remark}

Now consider $L$ and $L_0$ to be Legendrian in $(S^3, \xi_\textrm{st})$ and the surgeries to be contact $(\frac{1}{n})$-surgeries. We claim that the rotation number of $L_0$ in the surgered contact manifold $(M, \xi)$ with respect to $\widehat{\Sigma}$ is equal to
$$
\rot_M (L_0, \widehat{\Sigma}) = \rot_{S^3} (L_0) - \sum_{i=1}^k{a_i n_i \rot_{S^3} (L_i)}.
$$

In complete analogy to \cite{MR2557137}, \cite{MR3338830} and \cite{Conway2014} we have the following lemma.
\begin{lemma}
The homomorphism $H_1 (S^3 \setminus (L_0 \sqcup L)) \rightarrow H_1 (M \setminus L_0)$ induced by inclusion maps $\PD(\textrm{e}(\xi_\textrm{st}, L_0 \sqcup L ))$ to $\PD(\textrm{e}(\xi, L_0))$.
\end{lemma}
The proof is completely analogous to the ones in \cite{MR2557137, MR3338830}, where one uses the Legendrian rulings of the surgery tori induced by $(\frac{1}{n})$-surgery instead of $(\pm 1)$-surgery.

We thus have (cf.\ \cite{Conway2014})
\begin{align*}
\rot_M (L_0, \widehat{\Sigma}) = & \PD\big(\textrm{e}(\xi, L_0)\big) \bullet \widehat{\Sigma} \\
= & \PD\big(\textrm{e}(\xi_\textrm{st}, L_0 \sqcup L)\big) \bullet \Sigma \\
= & \left( \sum_{i=0}^k{ \rot_{S^3} (L_i) \mu_i} \right) \bullet \Sigma \\
= & \left( \sum_{i=0}^k{ \rot_{S^3} (L_i) \mu_i} \right) \bullet \left( \Sigma_0 + \sum_{j=1}^{k}{(-a_j n_j) \Sigma_j} \right) \\
= & \rot_{S^3} (L_0) - \sum_{i=1}^k{a_i n_i \rot_{S^3} (L_i)},
\end{align*}
which proves the theorem.
\end{proof}

If the contact surgeries are not unique, i.e.\ for contact surgery coefficients not of the form
$1/n$ (see Remark \ref{coefficients}), the rotation number is -- in contrast to the Thurston--Bennequin invariant -- not independent of the chosen contact structures on the solid tori, as the following example illustrates.

\begin{figure}[htbp] 
\centering
\def\svgwidth{0,65\columnwidth}
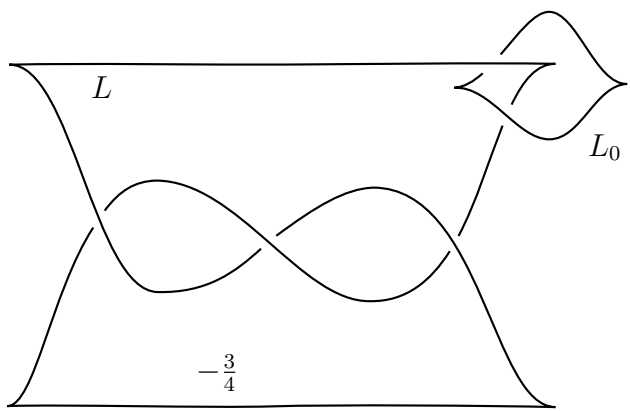
\caption{Non-unique contact surgery yielding a homology sphere}
\label{fig:Legtrefoil}
\end{figure}

\begin{example}
\label{ex:non-unique}
Consider the diagram depicted in Figure \ref{fig:Legtrefoil}, where $L$ is a Legendrian trefoil with contact surgery coefficient $3/4$ and $L_0$ a Legendrian unknot in its complement.
We have $\tb(L) = 1$, so the topological surgery coefficient is $\frac{1}{4}$. Thus, the surgered manifold $M$ is a homology sphere and the rotation number of $L_0$ independent of the choice of Seifert surface.
The contact surgery coefficient $-\frac{3}{4}$ has a continued fraction expansion $1 - 2 - \frac{1}{-4}$, which means that there are three distinct tight contact structures on the solid torus compatible with the surgery resulting in the contact manifolds which are shown in Figure \ref{fig:legtrefoilalgo} (see \cite{MR2056760}).
Topologically, these are the same, i.e.\ for all three diagrams we have
$$ Q =
\begin{pmatrix}
0 & 1 \\
1 & -2
\end{pmatrix}
\textrm{ and }
\mathbf{l} =
\begin{pmatrix}
1\\
1
\end{pmatrix},
$$
and hence $\mathbf{a} = (3, 1)$.
Furthermore, we have $\rot_{S^3} (L_0) = 0$, $\rot_{S^3} (L_1) = 0$ and $\rot_{S^3} (L_2) \in \{-2, 0, 2\}$.
This yields 
$$
\rot_M (L_0) = \rot_{S^3} (L_0) - 3 \rot_{S^3} (L_1) - \rot_{S^3} (L_2) \in \{-2, 0, 2\},
$$
depending on the chosen contact structure and orientations.

\begin{figure}[htbp] 
\centering
\def\svgwidth{0,9\columnwidth}
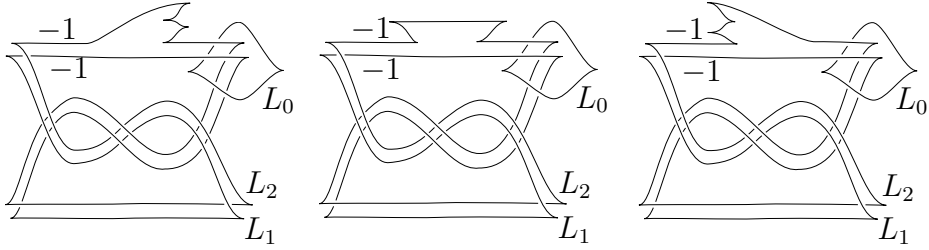
\caption{Three unique contact surgeries corresponding to Figure \ref{fig:Legtrefoil}}
\label{fig:legtrefoilalgo}
\end{figure}	

\end{example}

\begin{example}
We consider the case of $L$ a one component link with contact surgery coefficient $\pm\frac{1}{n}$, so the topological surgery coefficient is $\frac{(n\tb(L) \pm 1)}{n}$.
We then have $Q = p = n\tb(L) \pm 1$ and $L_0$ is nullhomologous in the surgered manifold if and only the linking number of $L_0$ and $L$ is divisible by $n\tb(L) \pm 1$, in which case $\mathbf{a}$ is the quotient $\frac{\lk (L_0,L)}{n\tb(L) \pm 1}$.
Then the rotation number of $L_0$ in the surgered manifold is
$$
\rot_M (L_0, \widehat{\Sigma}) = \rot_{S^3} (L_0) - { \frac{n \lk (L_0, L) }{n\tb(L) \pm 1} \rot_{S^3} (L)},
$$
and its Thurston--Bennequin invariant is
$$
\tb_{M} (L_0) = \tb_{S^3} (L_0) - { \frac{n \lk^2 (L_0, L)}{n\tb(L) \pm 1} }.
$$
Observe that if $n\tb(L) \pm 1$ is non-zero, the knot $L_0$ is rationally nullhomologous. Then the computed numbers represent the rational invariants (cf.\ Section \ref{section:rational}).
\end{example}

\begin{example}
Figure \ref{fig:stabilization} shows a stabilisation of a knot $L_0$. Topologically, the surgery along the meridian $L$ of $L_0$ again yields $S^3$, and there are two choices of tight contact structures on the solid torus compatible with the surgery. Here, in both cases, the resulting $S^3$ is tight. The topological knot type $L_0$ stays unchanged in $M$, but $L_0$ is either stabilised positively or negatively, depending on the particular choice of contact structure in the surgery torus.

The topological data in the diagrams with a unique choice is
$$ Q =
\begin{pmatrix}
0 & -1 \\
-1 & 3
\end{pmatrix}
\textrm{ and }
\mathbf{l} =
\begin{pmatrix}
1\\
1
\end{pmatrix},
$$
and hence $\mathbf{a} = (2, -1)$.
So the Thurston--Bennequin invariant of $L_0$ in $M$ is
$$
\tb_{M} (L_0) = \tb_{S^3} (L_0) - \big\langle \begin{pmatrix}2\\ -1\end{pmatrix} , \begin{pmatrix}1\\ 1\end{pmatrix} \big\rangle = \tb_{S^3} (L_0) - 1
$$
in both cases.
The rotation number of $L_1$ vanishes in both cases, the rotation number of $L_2$ is either $+1$ or $-1$.
We thus have
$$
\rot_M (L_0) = \rot_{S^3} (L_0) - \big\langle \begin{pmatrix}2\\ -1\end{pmatrix} , \begin{pmatrix}0\\ \pm 1\end{pmatrix} \big\rangle = \rot_{S^3} (L_0) \mp 1.
$$
In fact, one can show that it is a stabilised copy of $L_0$ (see \cite[Section~10]{Kegel2016}).

\begin{figure}[htbp] 
\centering
\def\svgwidth{0,9\columnwidth}
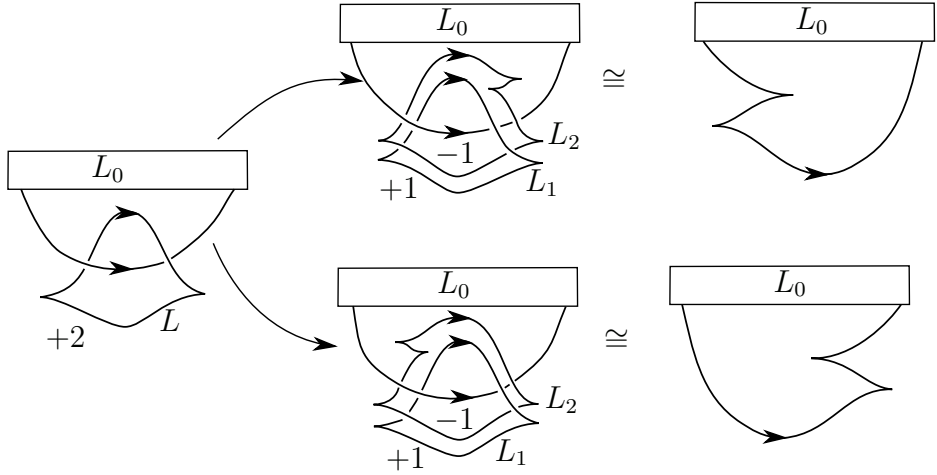
\caption{Stabilisation via surgery}
\label{fig:stabilization}
\end{figure}	

\end{example}

\section{The self-linking number of transverse knots}
\label{section:selfl}

Let $T$ be an oriented nullhomologous transverse knot in a contact manifold $(M, \xi)$ and let $\Sigma$ be a Seifert surface for $T$. The \emph{self-linking number} $\selfl (T, \Sigma)$ of $T$ is defined as the linking number of $T$ and $T'$ where $T'$ is obtained by pushing $T$ in the direction of a non-vanishing section of $\xi|_\Sigma$.

\begin{remark}
We consider transverse knots with arbitrary orientations. If the given orientation coincides with the orientation induced by the contact planes, we call the knot \emph{positively transverse}, and else \emph{negatively transverse}.
The self-linking number of a transverse knot is independent of its orientation and does only depend on the homology class of the chosen Seifert surface (cf.~Section~3.5.2 in \cite{Geiges2008}).
\end{remark}

\begin{corollary}
\label{cor:selfl}
Let $L = L_1 \sqcup \ldots \sqcup L_k$ be an oriented Legendrian link in $(S^3, \xi_\textrm{st})$ and $T_0$ an oriented transverse knot in its complement. Let $(M, \xi)$ be the contact manifold obtained from $S^3$ by contact $({1}/{n_i})$-surgeries ($n_i \in \Z$) along $L$.
Then $T_0$ is nullhomologous in $M$ if and only if there is an integral vector $\mathbf{a}$ solving $\mathbf{l} = Q\mathbf{a}$ as above, in which case its self-linking number in $(M, \xi)$ (with respect to the special Seifert class $\widehat{\Sigma}$ as before) is equal to
$$
\selfl_M (T_0, \widehat{\Sigma}) = \selfl_{S^3} (T_0) - \sum_{i=1}^k{a_i n_i \big(l_i \mp \rot_{S^3} (L_i)\big)},
$$
where the sign is $-$ when $T_0$ is positively transverse and $+$ when $T_0$ is negatively transverse.
\end{corollary}

\begin{remark}
An oriented transverse knot $T$ is either positively or negatively transverse. If we pick a Legendrian knot $L$ such that $T$ is a transverse push-off, we orient $L$ accordingly. Then the class of an oriented Seifert surface of $T$ is also the class of an oriented Seifert surface of $L$ and vice-versa. With these orientations, $T$ is a positive (negative) push-off of $L$ if $T$ is positively (negatively) transverse. In particular, the topological data used in the formula in Corollary \ref{cor:selfl} coincides for the two knots.
\end{remark}

\begin{proof}
Any transverse knot is a transverse push-off of a Legendrian knot (cf.\ the paragraph before Theorem 2.23 in \cite{MR2179261}), so it is enough to consider those.
Now for $L_0^\pm$ the positive or negative push-off of the Legendrian knot $L_0$ and $\Sigma$ a Seifert surface we have
$$
\selfl(L_0^\pm, [\Sigma]) = \tb(L_0) \mp \rot(L_0, [\Sigma])
$$
in any contact manifold (see Proposition 3.5.36 in \cite{Geiges2008}).
Hence,
\begin{align*}
\selfl_M (L_0^\pm, \widehat{\Sigma}) = & \tb_{M} (L_0) \mp \rot_{M} (L_0, \widehat{\Sigma}) \\
= & \Big( \tb_{S^3} (L_0) - \sum_{j=1}^{k}{a_j n_j l_{j}} \Big) \mp \Big( \rot_{S^3} (L_0) - \sum_{i=1}^k{a_i n_i \rot_{S^3} (L_i)} \Big) \\
= & \selfl_{S^3} (L_0) - \sum_{i=1}^k{a_i n_i \Big(l_i \mp \rot_{S^3} (L_i)\Big)}.
\end{align*}
\end{proof}

\begin{remark}
A front-projection that contains Legendrian as well as transverse knots has four possible types of crossings between a Legendrian and a transverse knot (see Figure \ref{fig:crossings}).
Depending on whether the transverse knot is positively or negatively transverse, two of the four types of crossings have a unique crossing behaviour determined by the contact condition, in the other cases both possibilities can occur.
\begin{figure}[htbp] 
\centering
\def\svgwidth{0,9\columnwidth}
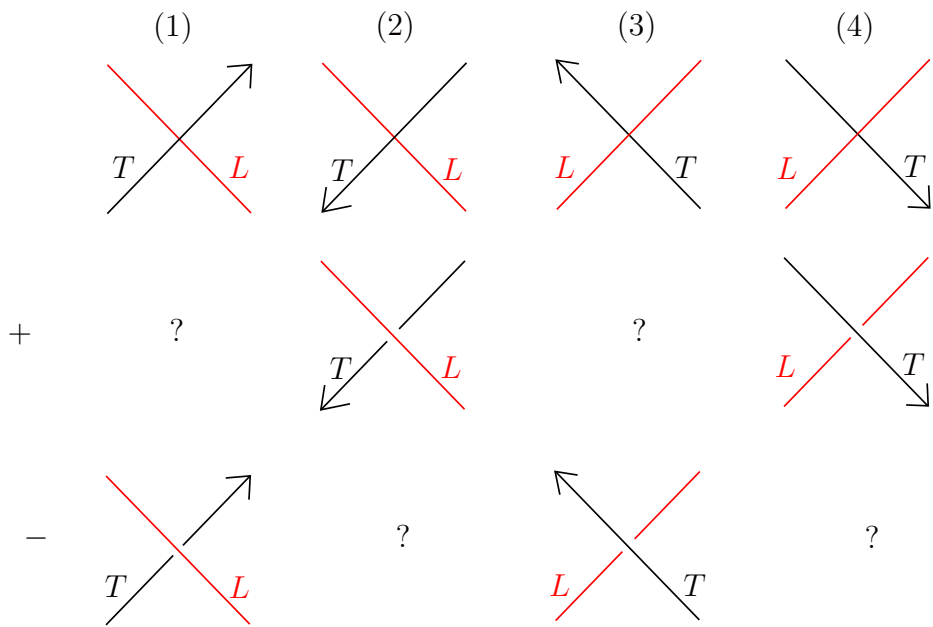
\caption{Crossings between Legendrian and transverse knots. The transverse knots in the middle row are positive, the ones in the bottom row negative.}
\label{fig:crossings}
\end{figure}
\end{remark}

\begin{example}
\begin{enumerate}
\item The left diagram in Figure~\ref{fig:slexample} shows a positive transverse knot $T_0$ in an overtwisted $3$-sphere $M$.
We have $\mathbf{l} = -1$, $Q = p = -1$ and thus $\mathbf{a} = 1$. The rotation number of $L$ is $1$, so we have
$$
\selfl_M(T_0) = \selfl_{S^3}(T_0) - a_1 q_1 (l_1 - \rot_{S^3}(L)) = -1 - (-1 - 1) = 1.
$$
Therefore, $T_0$ violates the Bennequin-inequality in $M$, i.e.\ the contact structure is indeed overtwisted.

Alternatively, we can consider a Legendrian unknot $L_0$ such that $T_0$ is its positive push-off, as shown on the right in Figure~\ref{fig:slexample}. Its Thurston--Bennequin invariant in $M$ is equal to $-1 + 1 = 0$ and its rotation number is $0-1=-1$, i.e.\ it bounds an overtwisted disc.

\begin{figure}[htbp] 
\centering
\def\svgwidth{0,9\columnwidth}
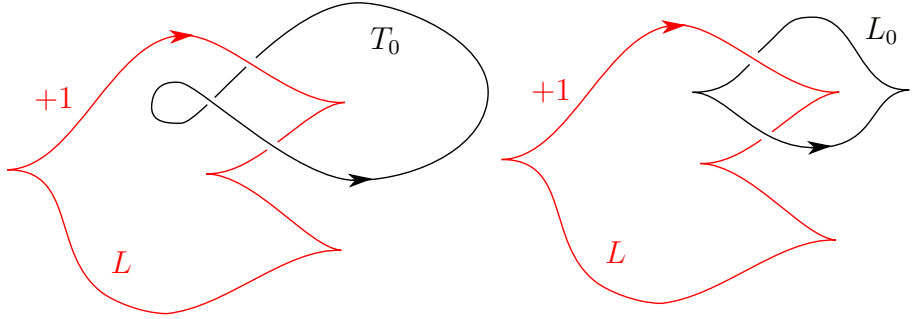
\caption{Computing the self-linking number}
\label{fig:slexample}
\end{figure}	

\item
We can also consider $T_0$ as a negative transverse knot by reversing its orientation. Then $\mathbf{l} = 1$, $Q = p = -1$ and  $\mathbf{a} = -1$, so 
$$
\selfl_M(T_0) = \selfl_{S^3}(T_0) - a_1 q_1 \big(l_1 - \rot_{S^3}(L)\big) = -1 - (1 + 1) = 1,
$$
as expected, since the self-linking number is independent of the chosen orientation.
We can again consider the corresponding Legendrian knot, which then has vanishing Thurston--Bennequin invariant and rotation number $1$. As $T_0$ is now its negative push-off, we also get
$$
\selfl_M(T_0) = \tb_M + \rot_M = 1.
$$
\end{enumerate}
\end{example}
\section{Rationally nullhomologous knots}
\label{section:rational}

The study of rationally nullhomologous knots in contact $3$-manifolds has been proposed in Baker-Grigsby~\cite{MR2552000}, Baker-Etnyre~\cite{MR2884030} and Geiges-Onaran~\cite{MR3338830}.
In this section we generalise Theorem \ref{thm:rot_surgery} to rationally nullhomologous Legendrian knots and Corollary \ref{cor:selfl} to rationally nullhomologous transverse knots.
Let $K$ be a knot in $M$. We call $K$ \emph{rationally nullhomologous} if its homology class is of finite order $d>0$ in $H_1(M)$, i.e.\ it vanishes in $H_1(M; \Q)$.
Let $\nu K$  be a tubular neighbourhood of $K$ and denote the meridian by $\mu \subset \partial\nu K$.

\begin{definition}
A \emph{Seifert framing} of an oriented rationally nullhomologous knot $K$ of order $d$ is a class $r \in H_1(\partial \nu K)$ such that
\begin{itemize}
\item $\mu \bullet r = d$,
\item $r = 0$ in $H_1(M \setminus \nu K)$.
\end{itemize}
A \emph{rational Seifert surface} for an oriented rationally nullhomologous knot $K$ is a surface with boundary in the complement of $K$ whose boundary represents a Seifert framing of $K$.
\end{definition}

It is obvious that every rationally nullhomologous knot has a Seifert framing. Moreover, the Seifert framing is unique (see \cite{tb_openbooks}).

\begin{definition}
The \emph{rational rotation number} of an oriented rationally nullhomologous Legendrian knot $K$ of order $d$ with respect to the rational Seifert surface $\Sigma$ is equal to
$$
\rot_\Q (K,\Sigma) = \frac{1}{d} \langle \textrm{e}(\xi, K), [\Sigma] \rangle = \frac{1}{d} \PD(\textrm{e}(\xi, K)) \bullet [\Sigma],
$$
where $\textrm{e}(\xi, K)$ is the relative Euler class of the contact structure $\xi$ relative to the knot $K$ and $[\Sigma]$ the relative homology class represented by the surface $\Sigma$ and the intersection is taken in $H_1(\partial \nu K)$.
\end{definition}

Let $L_0 \subset S^3\setminus L$ be an oriented knot in the complement of an oriented surgery link $L$. Using the notation from Section \ref{section:rotation}, we see that $L_0$ is rationally nullhomologous of order $d$ in $M = S^3_L(r)$  if and only if there is an integral solution $\mathbf{a}$ of the equation $d\mathbf{l} = Q\mathbf{a}$ and $d$ is the minimal natural number for which a solution exists (see \cite{Kegel2016}).

Now assume that $L$ and $L_0$ are Legendrian and $L_0$ is rationally nullhomologous of order $d$ in $M$ and fix Seifert surfaces $\Sigma_0, \ldots, \Sigma_k$ for $L_0, \ldots, L_k$ in $S^3$ such that intersections of surfaces and link components are transverse, as in the nullhomologous case.
Again following \cite{Conway2014}, we want to construct a class of the form
$$
\Sigma = d\Sigma_0 + \sum_{i=1}^{k}{k_i \Sigma_i}
$$
such that its image under the boundary homomorphism $\partial$ in the long exact sequence of the pair $(S^3 \setminus (L_0 \sqcup \nu L), \partial L_0 \sqcup \partial\nu L)$ is a linear combination of the surgery slopes $r_i$ and a rational longitude of $L_0$.
Setting $k_i = -a_i q_i$, we obtain
\begin{align*}
\partial \colon\thinspace \Sigma_j \longmapsto & \sum_{j=1}^{k}{d a_j q_j l_{0j}} \mu_0 + d \lambda_0 - \sum_{i=1}^{k}{a_i q_i \lambda_i} - \sum_{i=1}^k{d l_{0i} \mu_i} + \sum_{i=1}^{k}\sum_{j\neq i}{a_j q_j l_{ij} \mu_i} \\
= & \sum_{j=1}^{k}d {a_j q_j l_{j}} \mu_0 + d \lambda_0 - \sum_{i=1}^{k}{a_i q_i \lambda_i} - \sum_{i=1}^k{ dl_{i} \mu_i} + \sum_{i=1}^{k}{(d l_i - a_i  p_i )\mu_i} \\
= & \sum_{j=1}^{k}{d a_j q_j l_{j}} \mu_0 + d \lambda_0 - \sum_{i=1}^{k}{a_i q_i \lambda_i} - \sum_{i=1}^{k}{(a_i  p_i )\mu_i}.
\end{align*}

In complete analogy to the the nullhomologous case we then have
\begin{align*}
\rot_{\Q, M} (L_0, \widehat{\Sigma}) = & \frac{1}{d} \PD\big(\textrm{e}(\xi, L_0)\big) \bullet \widehat{\Sigma} \\
= & \rot_{S^3} (L_0) - \frac{1}{d} \sum_{i=1}^k{a_i n_i \rot_{S^3} (L_i)}.
\end{align*}

Thus, Theorem \ref{thm:rot_surgery} generalises as follows.

\begin{theorem}
In the situation of Theorem \ref{thm:rot_surgery} the knot $L_0$ is rationally nullhomologous of order $d$ in $M$ if and only if there is an integral vector $\mathbf{a}$ solving $d\mathbf{l} = Q\mathbf{a}$ as above with $d$ the minimal natural number for which a solution exists, in which case its rational rotation number in $(M, \xi)$ with respect to a special (rational) Seifert class $\widehat{\Sigma}$ is equal to
$$
\rot_{\Q, M} (L_0, \widehat{\Sigma}) = \rot_{S^3} (L_0) - \frac{1}{d}\sum_{i=1}^k{a_i n_i \rot_{S^3} (L_i)}.
$$
\end{theorem}

The definition of the self-linking number of a transverse knot generalises to the setting of rationally nullhomologous knots by choosing a rational Seifert surface.
Furthermore, the rational invariants of a Legendrian and the rational self-linking of a transverse push-off are, as in the nullhomologous case, related by
$$
\selfl_\Q(L_0^\pm, [\Sigma]) = \tb_\Q(L_0) \mp \rot_\Q(L_0, [\Sigma])
$$
(see Lemma 1.2 in \cite{MR2884030}). Hence, we have the following corollary.

\begin{corollary}
In the situation of Corollary \ref{cor:selfl} the knot $T_0$ is rationally nullhomologous of order $d$ in $M$ if and only if there is an integral vector $\mathbf{a}$ solving $d\mathbf{l} = Q\mathbf{a}$ as above, in which case its rational self-linking number in $(M, \xi)$ with respect to a special (rational) Seifert class $\widehat{\Sigma}$ is equal to
$$
\selfl_{\Q, M} (T_0, \widehat{\Sigma}) = \selfl_{S^3} (T_0) - \frac{1}{d} \sum_{i=1}^k{a_i n_i (l_i \mp \rot_{S^3} (L_i))}.
$$
\end{corollary}

\begin{remark}
Observe that the formulas for rationally nullhomologous knots coincide with the ones for nullhomologous knots presented in previous sections if one allows rational coefficients.
\end{remark}

\section{The $d_3$-invariant in surgery diagrams}
\label{section:d3}

The so-called $d_3$-invariant is a homotopical invariant of a
tangential $2$-plane field on a $3$-manifold, which is defined if the
Euler class (or first Chern class) of the $2$-plane field is torsion,
see \cite[Definition~11.3.3]{GS}.
Many contact structures can be distinguished by computing the $d_3$-invariants
of the underlying topological $2$-plane fields.
In \cite[Corollary~3.6]{MR2056760} Ding, Geiges and Stipsicz present a formula
to compute first the Euler class and then the $d_3$-invariant of a contact
structure given by a $(\pm1)$-contact surgery diagram building up on the work of
Gompf \cite{MR1668563}.
Both invariants are closely related to the rotation number of the surgery links.

By expressing an arbitrary $(1/n)$-contact surgery diagram as a
$(\pm1)$-contact surgery diagram and then using the result of
Ding--Geiges--Stipsicz we obtain a similar result for arbitrary $(1/n)$-contact
surgery diagrams, which often simplifies computations a lot.

First we recall some results from \cite{MR2056760}:
For $L = L_1 \sqcup \ldots \sqcup L_k$ an oriented Legendrian link in
$(S^3, \xi_\textrm{st})$ and $(M, \xi)$ the contact manifold
obtained from $S^3$ by contact $(\pm1)$-surgeries along $L$,
the Poicaré-dual of the Euler class is given by
\begin{equation*}
\operatorname{PD}\big(\textrm{e}(\xi)\big)=\sum_{i=1}^k \rot_i\mu_i\in H_1(M).
\end{equation*}
The meridians $\mu_i$ generate the first homology
$H_1(M)$ and the relations are given by $Q\mathbf\mu=0$.
Observe that he generalized linking matrix $Q$ coincides with the ordinary linking matrix, since we only have integer surgeries here.
Then $\textrm{e}(\xi)$ is torsion if and only if there exists a rational solution $\mathbf b\in\Q^k$ of $Q\mathbf b=\mathbf{rot}$.
If this is the case, then the $d_3$-invariant computes as
\begin{equation*}
\textrm{d}_3 = \frac{1}{4} \big(\langle \mathbf{b}, \mathbf{\rot} \rangle - 3 \sigma (Q)       - 2 k\big) - \frac{1}{2} + q,
\end{equation*}
where $\sigma (Q)$ denotes the signature of $Q$ (i.e.\ the number of positive eigenvalues minus the number of negative ones) and $q$ is the number of Legendrian knots in $L$ with $(+1)$-contact surgery coefficient.

With the help of these results we can now state and prove a corresponding theorem for arbitrary $(1/n)$-contact surgeries.


\begin{theorem}
Let $L = L_1 \sqcup \ldots \sqcup L_k$ be an oriented Legendrian link in
$(S^3, \xi_\textrm{st})$ and denote by $(M, \xi)$ the contact manifold obtained
from $S^3$ by contact $(\pm {1}/{n_i})$-surgeries along $L$ ($n_i \in \N$). 
\begin{enumerate}
	\item The Poicaré-dual of the Euler class is given by
\begin{equation*}
\operatorname{PD}\big(\textrm{e}(\xi)\big)=\sum_{i=1}^k n_i\rot_i\mu_i\in H_1(M).
\end{equation*}
The first homology group
$H_1(M)$ of $M$ is generated by the meridians $\mu_i$
and the relations are given by the generalized linking matrix $Q\mathbf\mu=0$.
\item The Euler class $\textrm{e}(\xi)$ is torsion if and only if there exists
a rational solution $\mathbf b\in\Q^k$ of $Q\mathbf b=\mathbf{rot}$.
In this case, the $d_3$-invariant computes as
\begin{equation*}
\textrm{d}_3 = \frac{1}{4} \left(\sum_{i=1}^k n_i b_i \rot_i  + (3-n_i) \operatorname{sign}_i\right)   - \frac{3}{4} \sigma (Q) - \frac{1}{2} ,
\end{equation*}
where $\operatorname{sign}_i$ denotes the sign of the contact surgery coefficient of $L_i$.
\end{enumerate}
\end{theorem}

\begin{remark}
In the proof we will show that all eigenvalues of $Q$ are real.
Therefore, it makes sense to speak of the signature, even if $Q$ is not symmetric.
\end{remark}

\begin{proof}
The replacement lemma of Ding and Geiges \cite[Proposition~8]{MR1823497} states that a contact $(\pm 1/n)$-surgery along a Legendrian knot $L$ is equivalent to $n$ contact $(\pm1)$-surgeries along Legendrian push-offs of $L$. Using this, we translate the contact $(\pm {1}/{n_i})$-surgeries along $L$ in contact $(\pm1)$-surgeries along a new Legendrian link $L'$ and compute the invariants there.

Denote by $L_i^j$ ($j=1,\ldots ,n_i$) the Legendrian push-offs of $L_i$ in the new Legendrian link $L'$.
Write $\mu_i$ for the meridian of $L_i$ ($i=1,\ldots,k$) and $\mu_i^j$ for the meridian of $L_i^j$ ($i=1,\ldots,k$, $j=1,\ldots ,n_i$).
We now have two surgery descriptions of the manifold $M$ -- one in terms of $L$
and one in terms of $L'$ -- and hence two presentations its first homology:
\begin{align*}
H_1(M)=&\langle\mu_i\vert Q\mathbf\mu=0\rangle \,\text{for the surgery presentation along $L$},\\
H_1(M)=&\langle\mu_i^j\vert Q'\mathbf\mu'=0\rangle \,\text{for the surgery presentation along $L'$}.
\end{align*}
An isomorphism between these two presentations is given by $\mu_i^j\mapsto \mu_i$ for
all $i,j$, and hence, as $\rot_i^j = \rot_i$,
$$
\operatorname{PD}\big(\textrm{e}(\xi)\big) = \sum_{i=1}^k\sum_{j=1}^{n_j} \rot_i^j\mu_i^j
\longmapsto
\sum_{i=1}^k n_i\rot_i\mu_i.
$$


The numbers $k$ and $q$ compute easily as
\begin{align*}
&k=\sum_{i=1}^k n_i,
&q=\sum_{i=1}^k \frac{1}{2}\left(1+\operatorname{sign}_i\right)n_i.
\end{align*}
For reasons of readability we will assume $k=3$ in the following. The general
case works exactly the same.
Write $\mathbf{1}_n$ for the vector $(1,\ldots,1)^T\in\Q^n$. 

Let $\mathbf b\in\Q^3$ a solution of $Q\mathbf b=\mathbf{rot}$, i.e\
\begin{align*}
Q\mathbf b=\begin{pmatrix} \pm1+n_1\tb_1& n_2l_{12} & n_3 l_{13} \\
n_1l_{12}&\pm1+n_2\tb_2 & n_3 l_{23} \\
n_1l_{13}& n_2l_{23} & \pm1+n_3\tb_3  \end{pmatrix}\begin{pmatrix} b_1 \\
b_2 \\
b_3 \end{pmatrix}=\begin{pmatrix} \rot_1 \\
\rot_2 \\
\rot_3 \end{pmatrix}=\mathbf{rot}
\end{align*}
Then for $\mathbf b':=(b_1,\ldots,b_1,b_2,\ldots,b_2,b_3,\ldots,b_3)^T\in\Q^{n_1+n_2+n_3}$ we have
\begin{align*}
Q'\mathbf b'=&\begin{pmatrix} \pm E_{n_1}+\tb_1 \mathbf{1}_{n_1}\mathbf{1}_{n_1}^T& l_{12} \mathbf{1}_{n_1}\mathbf{1}_{n_2}^T&  l_{13} \mathbf{1}_{n_2}\mathbf{1}_{n_3}^T\\
l_{12}\mathbf{1}_{n_2}\mathbf{1}_{n_1}^T&\pm E_{n_1}+\tb_2 \mathbf{1}_{n_2}\mathbf{1}_{n_2}^T &  l_{23} \mathbf{1}_{n_2}\mathbf{1}_{n_3}^T\\
l_{13}\mathbf{1}_{n_3}\mathbf{1}_{n_1}^T& l_{23}\mathbf{1}_{n_2}\mathbf{1}_{n_3}^T & \pm E_{n_3}+\tb_3 \mathbf{1}_{n_3}\mathbf{1}_{n_3}^T  \end{pmatrix}\mathbf b'\\
=&\begin{pmatrix} \rot_1 \mathbf{1}_{n_1} \\
\rot_2 \mathbf{1}_{n_2}\\
\rot_3 \mathbf{1}_{n_3} \end{pmatrix}=\mathbf{rot'}
\end{align*}
(conversely, every solution of $Q'\mathbf{b'} = \mathbf{rot'}$ is of this form and thus yields a solution of $Q\mathbf{b} = \mathbf{rot}$).
And therefore, 
\begin{equation*}
\langle \mathbf b',\mathbf{rot'}\rangle=\sum_{i=1}^3 n_i b_i \rot_i.
\end{equation*}

It remains to compute the signature $\sigma(Q')$ out of $\sigma(Q)$.
Let $\lambda$ be an eigenvalue of $Q$ with eigenvector $\mathbf v$.
Similar as above, 
one computes
\begin{equation*}
Q'\mathbf v'=\lambda\mathbf v'.
\end{equation*}
for $\mathbf v':=(v_1,\ldots,v_1,v_2,\ldots,v_2,v_3,\ldots,v_3)^T\in\Q^{n_1+n_2+n_3}$.
Thus, every eigenvalue of $Q$ is also an eigenvalue of $Q'$.
In particular, all eigenvalues of $Q$ are real.
Now we only have to find the other $\sum_{i=1}^3(n_i-1)$ eigenvectors of $Q'$.
To that end, consider the vector $\mathbf v_1\in \mathbf{1}_{n_1}^\bot$ and write $\mathbf v'_i:=(\mathbf v_1,0,\ldots,0,0,\ldots,0)^T\in\Q^{n_1+n_2+n_3}$. Then, as before, one computes
\begin{equation*}
Q'\mathbf v_1'=\operatorname{sign}_i\mathbf v_1',
\end{equation*}
and therefore
\begin{equation*}
\sigma(Q')=\sigma(Q)+\sum_{i=1}^k (n_i-1) \operatorname{sign}_i.
\end{equation*}
\end{proof}

\begin{example}
Consider a contact $(1/n)$-surgery ($n\in \Z$) along a Legendrian unknot with $\tb=-1$ and $\rot=0$. Then the Euler class is zero because the rotation number vanishes. Hence, the $\textrm{d}_3$-invariant is defined.
For $n = 1$, the signature of $Q$ vanishes.
If $n\neq 1$, the signature of $Q$ equals $-1$.
Thus we have
\begin{equation*}
\textrm{d}_3 =
\begin{cases}
\frac{n}{4} - \frac{1}{2}, & n < 1\\
0, & n = 1\\
1 - \frac{n}{4}, & n > 1\\
\end{cases}.
\end{equation*}
\end{example}

\section*{Acknowledgements}
The authors thank Christian Evers for helpful remarks on a draft version and Hendrik Herrmann for useful discussions about signatures and related topics. 

\end{document}